\documentclass{ifacconf}
\usepackage{graphicx}      % include this line if your document contains figures
\usepackage{natbib}        % required for bibliography
\usepackage{graphics} % for pdf, bitmapped graphics files
\usepackage{epsfig} % for postscript graphics files
\usepackage{amsmath} % assumes amsmath package installed
\usepackage{amssymb}  % assumes amsmath package installed
\usepackage{amsfonts}  % assumes amsmath package installed
\usepackage{subfig}
\usepackage{epstopdf}
\usepackage{enumerate}
\usepackage{graphicx} % grficos
\usepackage{color}

\newcommand*{\QEDB}{\hfill\ensuremath{\square}}%

\begin{document}
\begin{frontmatter}

\title{Fixed-Time Newton-Like Extremum Seeking} 

\thanks[footnoteinfo]{Research supported in part by CU Boulder - Autonomous Systems IRT Seed Grant No. 11005946, and NSF grants 1823983 $\&$ 1711373.}

\author[First]{J. I. Poveda} 
\author[Second]{M. Krsti\'{c}} 

\address[First]{Department of Electrical, Computer, and Energy Engineering, University of Colorado, Boulder, CO 80309 USA (e-mail: jorge.poveda@colorado.edu).}
\vspace{0.15cm}
\address[Second]{Department of Mechanical and Aerospace Engineering, University of California, San Diego, CA 92093 USA (e-mail: krstic@ucsd.edu).}

\begin{abstract}    % Abstract of not more than 250 words.
In this paper, we present a novel Newton-based extremum seeking controller for the solution of multivariable model-free optimization problems in static maps. Unlike existing asymptotic and fixed-time results in the literature, we present a scheme that achieves (practical) finite time convergence to a neighborhood of the optimal point, with a convergence time that is independent of the initial conditions and the Hessian of the cost function, and therefore can be arbitrarily assigned a priori by the designer with an appropriate choice of parameters in the algorithm. The extremum seeking dynamics exploit a class of fixed time convergence properties recently established in the literature for a family of Newton flows, as well as averaging results for perturbed dynamical systems that are not necessarily Lipschitz continuous. The proposed extremum seeking algorithm is model-free and does not require any explicit knowledge of the gradient and Hessian of the cost function. Instead, real-time optimization with fixed-time convergence is achieved by using real time measurements of the cost, which is perturbed by a suitable class of periodic excitation signals generated by a dynamic oscillator. Numerical examples illustrate the performance of the algorithm.
\end{abstract}


\end{frontmatter}
%===============================================================================

\section{Introduction}
In several applications it is of interest to recursively minimize a particular cost function whose mathematical form is unknown and which is only accessible via measurements. For these types of problems, extremum seeking (ES) algorithms have shown to be a powerful technique with provable stability, convergence, and robustness guarantees; see for instance \cite{KrsticBookESC,tan06Auto,DerivativesESC,Grushkovskaya2017,ThiagoKrstic,ESC_Suttner,Guay:03,Teel:01} and \cite{Poveda:16} for different extremum seeking architectures, applications, and theoretical results. Recently, there has been growing interest in improving the transient performance of ES algorithms by accelerating the rate of convergence. This has motivated the development in \cite[Sec. 6.1]{Poveda:16} of discontinuous algorithms, accelerated hybrid ES dynamics presented in \cite{zero_order_poveda_Lina}, as well as multivariable Newton-based schemes with an exponential rate of convergence that can be assigned a priory by the practitioner, see \cite{Newton}. Nevertheless, while these schemes exhibit better transient performance compared to the traditional gradient descent-based schemes, the convergence time is still highly dependent on the initial conditions of the algorithm. 

On the other hand, there has been a lot of recent efforts in designing control, estimation, and optimization algorithms with non-asymptotic convergence properties. These algorithms guarantee convergence to the desired target in a finite time that is independent of the initial conditions, see for instance  \cite{JaimeFixed_Time,finite_timeEngelTAC,Fixed_timeTAC,PralyFinite_Time,fixed_time,Fixed_timeTAC} and \cite{romeronips}. These results have opened the door to novel dynamical systems that are able to achieve fixed-time convergence using discontinuous or continuous vector fields; see for instance \cite{WeiRenFixed} and \cite{fixed_time}. Such results have also recently motivated a gradient-based ES algorithms with fixed-time convergence properties, presented in \cite{PovedaKrsticACC20}, where for the convex case the upper bound on the convergence time was shown to be independent of the initial conditions, but dependent on the smallest eigenvalue of the Hessian of the cost function, see \cite[Remark 1]{PovedaKrsticACC20}.

Motivated by this background, in this paper we present a novel Newton-based ES algorithm with fixed-time (practical) convergence properties, where the upper bound on the convergence time is now independent of the Hessian, and therefore can be arbitrarily assigned by the designer with a suitable choice of parameters in the algorithm. Our results are inspired by the previous architectures of Newton-based ES with asymptotic properties presented in \cite{Newton} and \cite{PowerES}, as well as by the recent model-based Newton-flows with fixed-time convergence properties introduced in \cite{fixed_time}.  In particular, we show that, under a suitable modification of the vector field used in standard multivariable Newton-based ES, fixed-time practical convergence can be achieved. Given that Newton-based ES dynamics estimate the inverse of the Hessian matrix of the cost function by using a Riccati differential equation that has multiple equilibria, our convergence results are local by nature. However, unlike existing results in the literature, the convergence time to an arbitrarily small neighborhood of the optimizer can be upper bounded by a positive number that is independent of the initial conditions and the Hessian of the cost function, and which can be prescribed a priori by the designer. This exhibits a clear advantage in comparison to traditional schemes whose convergence time depends heavily on the initial conditions and/or the Hessian of the cost. To our knowledge, the results of this paper correspond to the first Newton-based ES algorithm able to achieve (practical) fixed-time convergence with a fixed time that can be prescribed a priori.

The rest of this paper is organized as follows: In Section \ref{sec:preli} we present some preliminaries and definitions. In Section \ref{section3} we present the Newton-based extremum seeking controller considered in this paper, as well as the main convergence result. Section \ref{sec_apps} presents numerical simulations, and finally Section \ref{sec_conclusions} ends with the conclusions.
\section{PRELIMINARIES}
\label{sec:preli}
\subsection{Notation}
We denote the set of (non negative) real numbers by $(\mathbb{R}_{\geq 0})$ $\mathbb{R}$.
The set of (nonnegative) integers is denoted by $(\mathbb{Z}_{\geq 0})$ $\mathbb{Z}$. Given a compact
set $\mathcal{A} \subset \mathbb{R}^n$ and a vector $z \in \mathbb{R}^n$, we define $|z|_\mathcal{A} := \min_{y \in \mathcal{A}}|z-y|$, and we use $|\cdot|$ to denote the standard Euclidean norm. We denote by $
r\mathbb{B}$ the ball (in the Euclidean norm) of appropriate dimension centered around the origin and with radius $r>0$, and we use $\mathcal{A}+r\mathbb{B}$ to denote the set of points whose distance to $
\mathcal{A}$ is less or equal to $r$. For ease of notation, for two vectors $u,v \in \mathbb{R}^{n}$ we write $(u,v)$ for $(u^{T},v^{T})^{T}$. A function $\phi:\mathbb{R}^n\to\mathbb{R}$ is said to: a) be $C^k$ if its $k^{th}$ derivative is continuous; and b) be radially unbounded if $\phi(z)\to \infty$ whenever $|z|\to\infty$. We use $\mathbb{T}^n:=\mathbb{S}^1\times\ldots\times\mathbb{S}^1\subset\mathbb{R}^{2n}$ to denote the $n$-th cartesian product of unit circles centered around the origin, denoted as $\mathbb{S}^1:=\{u\in\mathbb{R}^2:u_1^2+u_2^2=1\}$. A function $\alpha:\mathbb{R}_{\geq0}\to\mathbb{R}$ is said to be of class $\mathcal{K}_{\infty}$ if it is zero at zero, continuous, strictly increasing, and unbounded. A function $\beta:\mathbb{R}_{\geq0}\times\mathbb{R}_{\geq0}\to\mathbb{R}$ is said to be of class $\mathcal{K}\mathcal{L}$ if it is nondecreasing in its first argument, non-increasing in its second argument, $\lim_{r\to0^+}\beta(r,s)=0$ for each $s\in\mathbb{R}_{\geq0}$, and $\lim_{s\to\infty}\beta(r,s)=0$ for each $r\in\mathbb{R}_{\geq0}$.

For the analysis of our algorithms, in this paper we will consider constrained dynamical systems of the form
\begin{equation}\label{diff_inclusion}
x\in C,~~~\dot{x}= F(x),
\end{equation}
where $F:\mathbb{R}^n\to\mathbb{R}^n$ is a continuous function, and $C\subset\mathbb{R}^n$ is a closed set. A continuously differentiable function $x:\text{dom}(x)\to\mathbb{R}^n$ is a solution to
 \eqref{diff_inclusion} if: (a) $x(0)\in C$, (b) $x(t)\in C$ for all $t\in\text{dom}(x)$, and (c) $\dot{x}(t)= F(x(t))$ for all $t\in\text{dom}(x)$. System \eqref{diff_inclusion} is said to render a compact set $\mathcal{A}\subset\mathbb{R}^n$ uniformly globally asymptotically stable (UGAS) if there exists a $\beta\in\mathcal{K}\mathcal{L}$ such that all solutions of \eqref{diff_inclusion} satisfy the bound 
 \begin{equation}
 |x(t)|_{\mathcal{A}}\leq \beta(|x(0)|_{\mathcal{A}},t),
\end{equation}
for all $t\in\text{dom}(x)$. When $\text{dom}(x)=[0,\infty)$ we say that the solution $x$ is complete.
\subsection{Problem Statement}
In this paper we are interested in solving the following unconstrained optimization problem
\begin{equation}\label{main_problem}
\min_{z\in\mathbb{R}^n}~~\phi(z),
\end{equation}
where $\phi:\mathbb{R}^n\to\mathbb{R}$ is an unknown cost function that is accessible only via measurements, and that satisfies the following assumption:

\begin{assum}\label{assumption_1}
The cost function $\phi$ is $C^3$, the Hessian $\nabla^2 \phi(z)$ is positive definite, and there exists $z^*\in\mathbb{R}^n$ such that $\nabla\phi(z^*)=0$. \QEDB
\end{assum}
In \cite{PovedaKrsticACC20}, it was shown that if $\nabla^2\phi(z)\geq \kappa I$, for some $\kappa>0$, then a gradient-based fixed-time extremum seeking algorithm can be used to solve problem \eqref{main_problem} in a finite time that can be taken to be independent of the initial conditions of the algorithm, but dependent on $\kappa$. In this paper, we are interested in removing this dependence by considering a Newton-like fixed-time extremum seeking (NFxTES) algorithm.

 %The radial unboundedness assumption on $|\nabla \phi(z)|^2$ will guarantee that some of the states of our algorithm can be initialized arbitrarily far away from $\mathcal{A}$. Therefore, this assumption can be relaxed if one is interested in establishing only local convergence results for all the states of the algorithm.

%Our goal is to design a Newton-based model-free optimization algorithm that guarantees convergence to $\mathcal{A}$ in a fixed time $T^*>0$ that is independent of the initial conditions.
%
\section{Newton-Based Extremum Seeking with Fixed-Time Convergence}
\label{section3}
In this section we present a novel \emph{model-free} Newton-like extremum seeking controller designed to solve problem \eqref{main_problem} in a fixed time.
\begin{figure}[t!]
 \begin{centering}
  \includegraphics[width=0.45\textwidth]{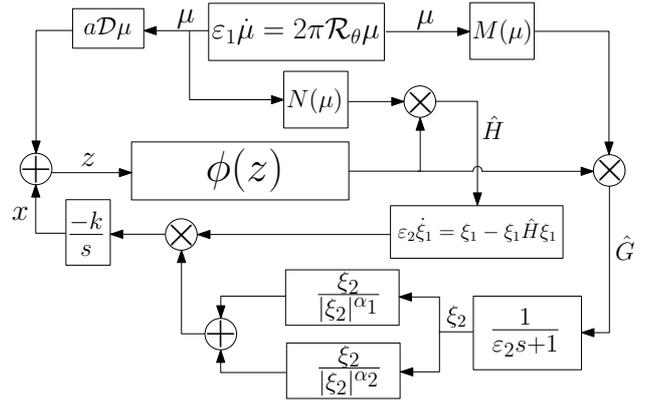} 
  \caption{Scheme of Newton-like Fixed-Time Extremum Seeking Dynamics for static maps.}
  \label{FigScheme}
 \end{centering}
\end{figure}
\subsection{Extremum Seeking Architecture}
Consider the Newton-based fixed-time extremum seeking (NFxTES) algorithm shown in Figure \ref{FigScheme}, and characterized by the equations:
\begin{equation}\label{ES_dynamics}
\dot{x}=F_x(\xi_1,\xi_2):=-\xi_1\xi_2\left(\frac{k}{|\xi_2|^{\alpha_1}}+\frac{k}{|\xi_2|^{\alpha_2}}\right),\\
\end{equation}
where the parameters $\alpha_1$ and $\alpha_2$ are defined as
\begin{equation}\label{alphas}
\alpha_1:=\frac{q_1-2}{q_1-1},~~~~\alpha_2:=\frac{q_2-2}{q_2-1},
\end{equation}
and where $(q_1,q_2,k)\in\mathbb{R}_{>0}^3$ are tunable parameters that will have an important role in the convergence properties of the algorithm.

The states $\xi_1$ and $\xi_2$ are generated by the dynamics
\begin{equation}\label{estimator_Hessian}
\dot{\xi}_1=F_{\xi_1}(\xi,x,\mu):=\frac{1}{\varepsilon_2}\Big(\xi_1-\xi_1\phi(z)N(\mu) \xi_1\Big),
\end{equation}
and
\begin{equation}\label{estimator_Gradient}
\dot{\xi}_2=F_{\xi_2}(\xi,x,\mu):=\frac{1}{\varepsilon_2}\Big(-\xi_2 +\phi(z)M(\mu)\Big),
\end{equation}
where $\varepsilon_2\in\mathbb{R}_{>0}$ is a tunable parameter selected sufficiently small to guarantee that \eqref{estimator_Hessian} and \eqref{estimator_Gradient} evolve in a faster time scale compared to the dynamics \eqref{ES_dynamics}. As shown in the next section, the state $\xi_1$ can be seen as an estimation of the inverse of the Hessian matrix $\nabla^2\phi(x)$, while the state $\xi_2$ can be seen as an estimation of the gradient $\nabla \phi(x)$.

In system \eqref{estimator_Hessian}-\eqref{estimator_Gradient}, the state $\mu$ is generated by a constrained linear oscillator evolving on $\mathbb{T}^n$, given by:
\begin{equation}\label{equation3}
\dot{\mu}=\frac{2\pi}{\varepsilon_1}\mathcal{R}_{\theta}\mu,~~~~~~\mu\in\mathbb{T}^n,
\end{equation}
where the parameter $\varepsilon_1>0$ is selected sufficiently small to guarantee that the dynamics \eqref{equation3} operate in a faster time scale compared to the dynamics \eqref{ES_dynamics}, \eqref{estimator_Hessian}, and \eqref{estimator_Gradient}, i.e., $\varepsilon_1\ll\varepsilon_2$. In order to guarantee that \eqref{equation3} is a dynamic oscillator, we define the matrix $\mathcal{R}_{\theta}\in\mathbb{R}^{2n\times2n}$ as a block diagonal matrix with blocks $\mathbf{R}_i\in\mathbb{R}^{2\times 2}$ defined as
\begin{equation*}
\mathbf{R}_i:=\left[\begin{array}{cc}
0 & \theta_i\\
-\theta_i & 0
\end{array}\right],~~~\theta_i\in\mathbb{R}_{>0}.
\end{equation*}
To link the states $(x,\mu)$, generated by the dynamics \eqref{ES_dynamics} and \eqref{equation3}, to the argument of the cost function $\phi(z)$, we set the argument $z$ via the feedback law
\begin{equation}\label{input}
z:=x+a\tilde{\mu},
\end{equation}
where $a\in\mathbb{R}_{>0}$ is a tunable parameter, and where $\tilde{\mu}\in\mathbb{R}^n$ is the vector that contains the odd entries of  $\mu\in\mathbb{R}^{2n}$, i.e.,
\begin{equation}\label{tildemu}
\tilde{\mu}:=[\mu_1,\mu_3,\mu_5,\ldots,\mu_{2n-1}]^\top.
\end{equation}
Note that $\tilde{\mu}$ can be written as $\tilde{\mu}=\mathcal{D}\mu$, with $\mathcal{D}\in\mathbb{R}^{n\times 2n}$ being a suitable matrix with entries $\mathcal{D}_{i,j}\in\{0,1\}$.

The gradient and Hessian estimators \eqref{estimator_Hessian}-\eqref{estimator_Gradient} depend on the excitation functions $N:\mathbb{R}^{2n}\to\mathbb{R}^{n\times n}$ and $M:\mathbb{R}^{2n}\to\mathbb{R}^n$, which are defined as
\begin{equation}\label{excitation_signals}
M(\mu):=\frac{2}{a}\mathcal{D}\mu,~~~N(\mu):=\left[\begin{array}{cccc}
N_{11} & N_{12} & \ldots & N_{1n}\\
N_{21} & N_{22} & \ldots & N_{2n}\\
\vdots  & \vdots & \vdots & \vdots\\
N_{n1} & N_{n2} & \ldots & N_{nn}
\end{array}\right],
\end{equation}
where the entries $N_{ij}$ satisfy $N_{ij}=N_{ji}$, as well as the conditions:
\begin{subequations}
\begin{align}
N_{ij}&=\frac{16}{a^2}\left(\tilde{\mu}_i^2-\frac{1}{2}\right),~~~~~\forall~i=j,\label{signal_Nii}\\
N_{ij}&=\frac{4}{a^2}\tilde{\mu}_i\tilde{\mu}_j,~~~~~~~~~~~~~~\forall~i\neq j\label{signal_Nij}.
\end{align}
\end{subequations}
As shown in the next section, under suitable choices of the parameters $(a,\varepsilon_1,\theta_i)$, the excitation signals \eqref{excitation_signals} guarantee an estimation of order $O(a)$ of the Hessian $\nabla^2\phi(x)$ and the gradient $\nabla\phi(x)$ via the dynamics \eqref{estimator_Hessian} and \eqref{estimator_Gradient}. Since $a$ is a tunable parameter, the estimation error can be made arbitrarily small on compact sets.
\begin{rem}
The structure of the NFxTES is similar to the multi-variable Newton-based extremum seeking controller considered in \cite{Newton}. However, the NFxTES dynamics have two main differences: a) under suitable choices of the parameters $(q_1,q_2)$, the dynamics of $x$ are continuous but not Lipschitz continuous, and they aim to approximate a Newton-based flow with fixed-time convergence properties instead of the Newton-based flow $\dot{x}=-\nabla^2\phi(x)^{-1}\nabla\phi(x)$ considered in \cite{Newton}; b) the dither signals $N(\cdot)$, $M(\cdot)$ and $\tilde{\mu}(\cdot)$ are generated by the linear oscillator \eqref{equation3}, which evolves on the torus $\mathbb{T}^n$. This allows us to analyze the NFxTES as a \emph{time-invariant} dynamical system. \QEDB
\end{rem}
\subsection{Main Result}
In this section, we establish the main convergence properties of the NFxTES. In order to do this, we define the \textsl{admissible} parameters $(q_1,q_2,k)$ as any tuple that satisfies the following relationships:
\begin{equation}\label{relationships}
q_1\in(2,\infty),~~q_2\in(1,2),~~k\in(0,\infty)
\end{equation}
As shown in \cite[Lemma 6]{fixed_time}, an admissible tuple of parameters guarantees that $1-\alpha_1>0$ and $1-\alpha_2>0$. Therefore, we have that
\begin{align*}
\left|\lim_{\xi_2\to0}\frac{\xi_1\xi_2}{|\xi_2|^{\alpha_1}}\right|\leq \left| \xi_1\right|\lim_{\xi_2\to0} |\xi_2|^{1-\alpha_1}=0.
\end{align*}
Using the same procedure for $\left|\frac{\xi_1\xi_2}{|\xi_2|^{\alpha_2}}\right|$, one obtains that $\left|\lim_{\xi_2\to0}\frac{\xi_1\xi_2}{|\xi_2|^{\alpha_2}}\right|=0$, which establishes continuity of the vector field \eqref{ES_dynamics} at points satisfying $\xi_2=0$. Continuity at points $\xi_2\neq0$ follows trivially by the structure of the dynamics. Based on this, we consider the following assumption on the parameters of the algorithm:
\begin{assum}\label{assumption_frequencies}
The tuple of parameters $(q_1,q_2,k)$ is admissible; and for each $i\in\{1,2,\ldots,n\}$ the parameter $\theta_i$ is a positive rational number, and $\theta_i\neq \theta_j$ for all $i\neq j$.  \QEDB
\end{assum}
\begin{rem}
The second condition of Assumption \ref{assumption_frequencies} is standard in extremum seeking control, see for instance \cite{DerivativesESC}, \cite{Newton}, and \cite{Poveda:16}, and it facilitates the application of averaging theory to analyze the qualitative behavior of the extremum seeking dynamics.~~~\QEDB
\end{rem}

For each admissible tuple of parameters $(q_1,q_2,k)$, we define the value
\begin{equation}\label{upper_bound}
T^*:=\frac{1}{k}\left[\left(\frac{2^{0.5\alpha_1}}{\alpha_1}\right)-\left(\frac{2^{0.5\alpha_2}}{\alpha_2}\right)\right],
\end{equation}
where $(\alpha_1,\alpha_2)$ are defined as in \eqref{alphas}. Note that admissible parameters always guarantee that the term inside the brackets is positive. Thus, for each desired $T^*>0$ one can always satisfy equation \eqref{upper_bound} by choosing admissible parameters $(q_1,q_2,k)$ with 
\begin{equation}\label{choiceC}
k=\frac{1}{T^*}\left[\left(\frac{2^{0.5\alpha_1}}{\alpha_1}\right)-\left(\frac{2^{0.5\alpha_2}}{\alpha_2}\right)\right].
\end{equation}
Using the definition of $T^*$ in \eqref{upper_bound}, and the NFxTES dynamics \eqref{ES_dynamics}, \eqref{estimator_Hessian}, \eqref{estimator_Gradient}, and \eqref{equation3}, we are now ready to state the main result of the paper. The complete proof is presented in \cite{fixed_timeES_arxiv}.
\begin{thm}\label{main_theorem}
Consider the NFxTES dynamics and suppose that Assumptions \ref{assumption_1} and \ref{assumption_frequencies} hold. Then, for admissible parameters $(q_1,q_2,k)$ and each $\nu>0$ there exists $\varepsilon_2^*>0$ such that for each $\varepsilon_2\in(0,\varepsilon_2^*)$ there exists $a^*>0$ such that for each $a\in(0,a^*)$ there exists $\varepsilon_1^*>0$ such that for each $\varepsilon_1\in(0,\varepsilon_1^*)$ there exists a neighborhood $\mathcal{N}$ of $p^*:=(z^*, \nabla^2 \phi(z^*)^{-1},\nabla\phi(z^*))$ such that every solution with $(x(0),\xi(0),\mu(0))\in \mathcal{N}\times\mathbb{T}^n$ exists for all time, and satisfies 
\begin{equation}\label{finite_convergence_bound}
|x(t)-z^*|\leq \nu,~~~\forall~t\geq T^*,
\end{equation}
where $T^*$ is given by \eqref{upper_bound}. \QEDB
\end{thm}

While the convergence result of Theorem \ref{main_theorem} is local with respect to the initial conditions, practical with respect to the parameters $(\varepsilon_1,\varepsilon_2,a)$ and the neighborhood $\{z^*\}+\nu\mathbb{B}$, and finite with respect to time, there is a key difference with respect to previous results in the literature of extremum seeking: the upper bound $T^*$ on the convergence time is independent of the initial conditions and the Hessian of the cost function, and therefore can be arbitrarily assigned a priori by the designer using an appropriate choice of admissible parameters $(q_1,q_2,k)$, such as in \eqref{choiceC}. This represents a clear advantage even with respect to \emph{gradient}-based fixed-time ES dynamics, e.g., \cite{PovedaKrsticACC20}, where the value of $T^*$ depends on the smallest eigenvalue of the Hessian matrix $\nabla^2\phi(z)$.

\textbf{Sketch of the Proof of Theorem 5:} We start with the following auxiliary Lemma:
\begin{lem}\label{integrals2}
Suppose that Assumption \ref{assumption_frequencies} holds and consider the linear oscillator \eqref{equation3} with $\varepsilon_1=1$ and the signal $\tilde{\mu}$ defined in \eqref{tildemu}. Then, there exists a $T>0$ such that for all solutions $\mu$ of \eqref{equation3} the following integrals hold:
\begin{subequations}\label{integrals}
\begin{align}
&\frac{1}{T}\int_{0}^T\tilde{\mu}_i(s)ds=0,~~\forall~~i\in\{1,2,\ldots,n\}.\\
&\frac{1}{T}\int_{0}^T\tilde{\mu}_i(s)^2ds=\frac{1}{2},\\
&\frac{1}{T}\int_{0}^T\tilde{\mu}_{i}(s)\tilde{\mu}_j(s)ds=0,~~\forall~~i\neq j,\\
&\frac{1}{T}\int_{0}^T\tilde{\mu}_{i}(s)^2N_{ii}(s)ds=\frac{1}{8},~~\forall~~i\in\{1,2,\ldots,n\},\\
&\frac{1}{T}\int_{0}^T\tilde{\mu}_{i}(s)^2N_{jj}(s)ds=0,~~\forall~~i\neq j,\\
&\frac{1}{T}\int_{0}^T\tilde{\mu}_{i}(s)^2N_{ij}(s)ds=0,~~\forall~~i\neq j,\\
&\frac{1}{T}\int_{0}^T\tilde{\mu}_{i}(s)N_{ij}(s)ds=0,~~~\forall~~i\neq j,\\
&\frac{1}{T}\int_{0}^T\tilde{\mu}_{i}(s)N_{ii}(s)ds=0,~~\forall~~i\in\{1,2,\ldots,n\},\\
&\frac{1}{T}\int_{0}^T\tilde{\mu}_{i}(s)N_{jj}(s)ds=0,~~\forall~~i\neq j,\\
&\frac{1}{T}\int_{0}^T\tilde{\mu}_{i}(s)\tilde{\mu}_{j}(s)N_{ij}(s)ds=\frac{1}{4},~~\forall~~i\neq j,\\
&\frac{1}{T}\int_{0}^T\tilde{\mu}_{i}(s)\tilde{\mu}_{j}(s)N_{ii}(s)ds=0,~~\forall~~i\neq j.
\end{align}
\end{subequations}
  \QEDB
\end{lem}
\textbf{Proof:} The proof of Lemma \ref{integrals2} follows by taking $T=2\pi\times LCM\{\frac{1}{\kappa_1},\frac{1}{\kappa_2},\ldots,\frac{1}{\kappa_n}\}$, where LCM stands for the least common multiplier, and noticing that every solution $\tilde{\mu}_i$ generated by the oscillator \eqref{equation3} with $\varepsilon_1=1$ is of the form
\begin{equation}\label{solutions_oscillator}
\tilde{\mu}_i(t)=\mu_{i,1}(0)\cos(2\pi \kappa_it)+\mu_{i,2}(0)\sin(2\pi \kappa_it),
\end{equation}
with $\mu_{i,1}(0)^2+\mu_{i,2}(0)^2=1$. Using the expression \eqref{solutions_oscillator} all the integrals of \eqref{integrals} follow by direct computation. \null \hfill \null $\blacksquare$

The properties of Lemma \ref{integrals2} allow us to obtain real-time gradient and Hessian estimations in \eqref{estimator_Hessian} and \eqref{estimator_Gradient} -on average- by using only measurements of the cost function $\phi(z)$. In particular, by taking a Taylor expansion of $\phi(x+a\tilde{\mu})$ around $x$ for small values of $a$, we obtain:
\begin{equation*}
\phi(x+a\tilde{\mu})=\phi(x)+a\tilde{\mu}^\top \nabla \phi(x)+\frac{a^2}{2}\tilde{\mu}^\top\nabla^2 \phi(x)\tilde{\mu}+O(a^3).
\end{equation*}
Using the definitions of $M(\cdot)$ and $N(\cdot)$, and the integrals \eqref{integrals}, we obtain the following average functions, where the average is taken with respect to the solutions of the oscillator \eqref{equation3}, i.e., keeping $x$ constant:
\begin{equation}\label{average1}
\frac{1}{T}\int_{0}^{T}\phi(x+a\tilde{\mu}(s))M(\mu(s))ds=\nabla\phi(x)+O(a),
\end{equation}
and
\begin{equation}\label{average2}
\frac{1}{T}\int_{0}^{T}\phi(x+a\tilde{\mu}(s))N(\mu(s))ds=\nabla^2\phi(x)+O(a),
\end{equation}
which correspond to $O(a)$-perturbed versions of the gradient and the Hessian of the cost function $\phi$, with a perturbation that shrinks as $a\to0^+$. Therefore, by computing the average of the dynamics \eqref{ES_dynamics}, \eqref{estimator_Hessian}, and \eqref{estimator_Gradient} along the solutions of \eqref{equation3}, and by neglecting the $O(a)$ disturbance of \eqref{average1} and \eqref{average2}, we can obtain the following average dynamics:
\begin{subequations}\label{ES_dynamics2}
\begin{align}
\dot{x}&=F^A_x(\xi_1,\xi_2):=-\xi_1\xi_2\left(\frac{k}{\|\xi_2\|^{\alpha_1}}+\frac{k}{\|\xi_2\|^{\alpha_2}}\right),\label{opti_dynamics2}\\
\dot{\xi}_1&=F^A_{\xi_1}(\xi,x):=\frac{1}{\varepsilon_2}\Big(\xi_1-\xi_1\nabla^2\phi(x)\xi_1\Big),\label{hessian_dynamics}\\
\dot{\xi}_2&=F^A_{\xi_2}(\xi,x):=\frac{1}{\varepsilon_2}\Big(-\xi_2 +\nabla \phi(x)\Big).\label{gradient_dynamics}
\end{align}
\end{subequations}
When $\varepsilon_2$ is sufficiently small, system \eqref{ES_dynamics2} is a singularly perturbed system. The boundary layer dynamics of this system can be obtained by setting $\dot{x}=0$ and $\varepsilon_2=1$:
\begin{subequations}\label{ES_dynamics3}
\begin{align}
\dot{\xi}_1&=\xi_1-\xi_1\nabla^2\phi(x)\xi_1,\label{hessian_dynamics3}\\
\dot{\xi}_2&=-\xi_2 +\nabla \phi(x).\label{gradient_dynamics3}
\end{align}
\end{subequations}
For fixed values of $x$, the stability properties of system \eqref{ES_dynamics3} can be analyzed as follows: denote by $r_x:=\nabla \phi(x)$ and $H_x:=\nabla^2\phi(x)$, and consider the errors $\tilde{\xi}_1=\xi_1-H_x^{-1}$ and $\tilde{\xi}_2=\xi_2-r_x$. The error boundary layer dynamics are then given by the following decoupled equations:
\begin{subequations}
\begin{align}
\dot{\tilde{\xi}}_1&=-\tilde{\xi}_1H_x(\tilde{\xi}_1+H_x^{-1})\label{filter1}\\
\dot{\tilde{\xi}}_2&=-\tilde{\xi}_2.\label{filter2}
\end{align}
\end{subequations}
System \eqref{filter2} has the origin globally exponentially stable, and system \eqref{filter1} has the origin locally exponentially stable since its linearization around the origin has the Jacobian $-I$, see \cite[pp. 1761]{Newton}. Therefore, for each fixed $x$, the boundary layer dynamics render the quasi-steady state mapping
\begin{equation*}
 \xi^*=M(x):=(\nabla^2 \phi(x)^{-1},\nabla \phi(x)),
\end{equation*}
locally exponentially stable, uniformly on $x$.

%\begin{equation}\label{linear_system}
%\left[\begin{array}{c}
%\dot{\tilde{\xi}}_1\\
%\dot{\tilde{\xi}}_2
%\end{array}\right]=\left[\begin{array}{cc}
%-H &~~H^{-1}\\
%~0 & -I
%\end{array}\right]\left[\begin{array}{c}
%\tilde{\xi}_1\\
%\tilde{\xi}_2
%\end{array}\right],
%\end{equation}
%%
%Since the determinant of the matrix of \eqref{linear_system} is given by $\text{det}(-H)\text{det}(-I)$, by Assumption \ref{assumption_1} we conclude that system is exponentially stable, uniformly on $x$.

The reduced dynamics associated to system \eqref{ES_dynamics2} are obtained by substituting $\xi^*$ in \eqref{opti_dynamics2}:
\begin{align}\label{slow_dynamics}
\dot{x}=-k\nabla^2\phi(x)^{-1}\left(\frac{\nabla\phi(x)}{\|\nabla\phi(x)\|^{\alpha_1}}+\frac{\nabla\phi(x)}{\|\nabla\phi(x)\|^{\alpha_2}}\right).
\end{align}
Following the ideas of \cite{fixed_time}, system \eqref{slow_dynamics} can be analyzed using the Lyapunov function
\begin{equation*}
V(x)=\frac{1}{2}|\nabla \phi(x)|^2,
\end{equation*}
which under Assumption \ref{assumption_1} is positive definite with respect to the point $\{z^*\}$, and also has bounded level sets. The derivative of $V$ along the solutions of \eqref{slow_dynamics} is given by
\begin{equation}\label{Lyapunov}
\dot{V}(x)=-c_1 2^{\frac{\tilde{\alpha}_1}{2}}V(x)^{\frac{\tilde{\alpha}_1}{2}}-c_2 2^{\frac{\tilde{\alpha}_2}{2}}V(x)^{\frac{\tilde{\alpha}_2}{2}},
\end{equation}
which, by \cite[Lemma 1]{Fixed_timeTAC}, implies that the point $\{z^*\}$ is fixed-time stable. Thus, all solutions of \eqref{slow_dynamics} satisfy a bound of the form
\begin{equation*}
|x(t)-z^*|\leq \beta_x(|x(0)-z^*|,t),~~\forall~t\geq0,
\end{equation*}
where $\beta_x\in\mathcal{K}\mathcal{L}_{\mathcal{T}}$, i.e., $\beta$ is a $\mathcal{K}\mathcal{L}$ function that also satisfies the condition \cite[Sec. 2.1]{RiosTeel18} $\beta_x(r,t)=0$, for all $t\geq T^*$ and $r\geq0$, where $T^*$ is given by \eqref{upper_bound}. From here, the main result of the theorem follows by applying the averaging results of \cite[Thm. 7]{zero_order_poveda_Lina} that preserves the $\mathcal{K}\mathcal{L}$ bound of the slow dynamics \eqref{slow_dynamics} for the evolution of the state $x$ in the original dynamics \eqref{ES_dynamics}. \null \hfill \null $\blacksquare$

\begin{rem}
The local nature of Theorem \ref{main_theorem} is due to the existence of multiple equilibria in the dynamics
\begin{equation}\label{estimator_Hessian2}
\dot{\xi}_1=\Big(\xi_1-\xi_1\nabla^2\phi(z) \xi_1\Big),
\end{equation}
which corresponds to the average dynamics of \eqref{estimator_Hessian}. Similar local results emerge in Newton-based ESCs with asymptotic convergence properties; see for instance \cite{Newton}. While it is possible to design Newton-based ESCs with semi-global practical asymptotic stability results by computing the vector $\nabla^2\phi(x)^{-1}\nabla\phi(x)$ using the singular perturbation approach presented in \cite[Sec. 3]{NewtonSemiglobal}, said approach cannot be used in this case since it will generate dynamics \eqref{ES_dynamics} with discontinuous vector fields that are not locally bounded. % Nevertheless, even though Theorem \ref{main_theorem} is of local nature, as shown in Section \ref{sec_apps}, the NFxTES dynamics can achieve convergence to the solution of \eqref{main_problem} from initial conditions $x(0)$ that are not necessarily close to the optimizer $z^*$, provided the radial unboundedness condition of Assumption \ref{assumption_1} holds.  
\QEDB
\end{rem}

%%%%%%%%%%%%%%%%%%%%%%%%%%%%%%%%%%%%%%%%%%%
\section{NUMERICAL RESULTS}
\label{sec_apps}
\begin{figure}[t!]
 \begin{centering}
  \includegraphics[width=0.45\textwidth]{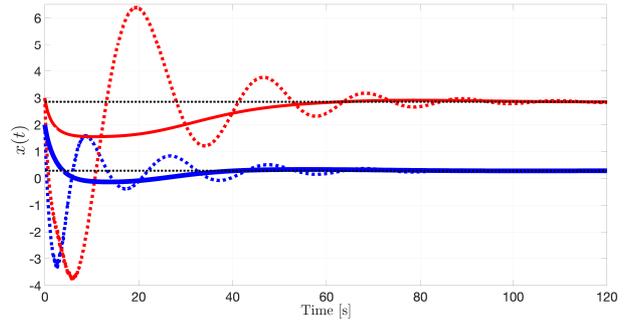}
  \caption{Evolution in time of the state $x$. Blue line corresponds to $x_1$, and red line corresponds to $x_2$. The dotted lines correspond to the trajectories generated by the traditional Newton-based ESC of \cite{Newton}. The solid lines correspond to the trajectories generated by the NFxTES.}
  \label{Fig1}
 \end{centering}
\end{figure}

\begin{figure}[t!]
 \begin{centering}
  \includegraphics[width=0.24\textwidth]{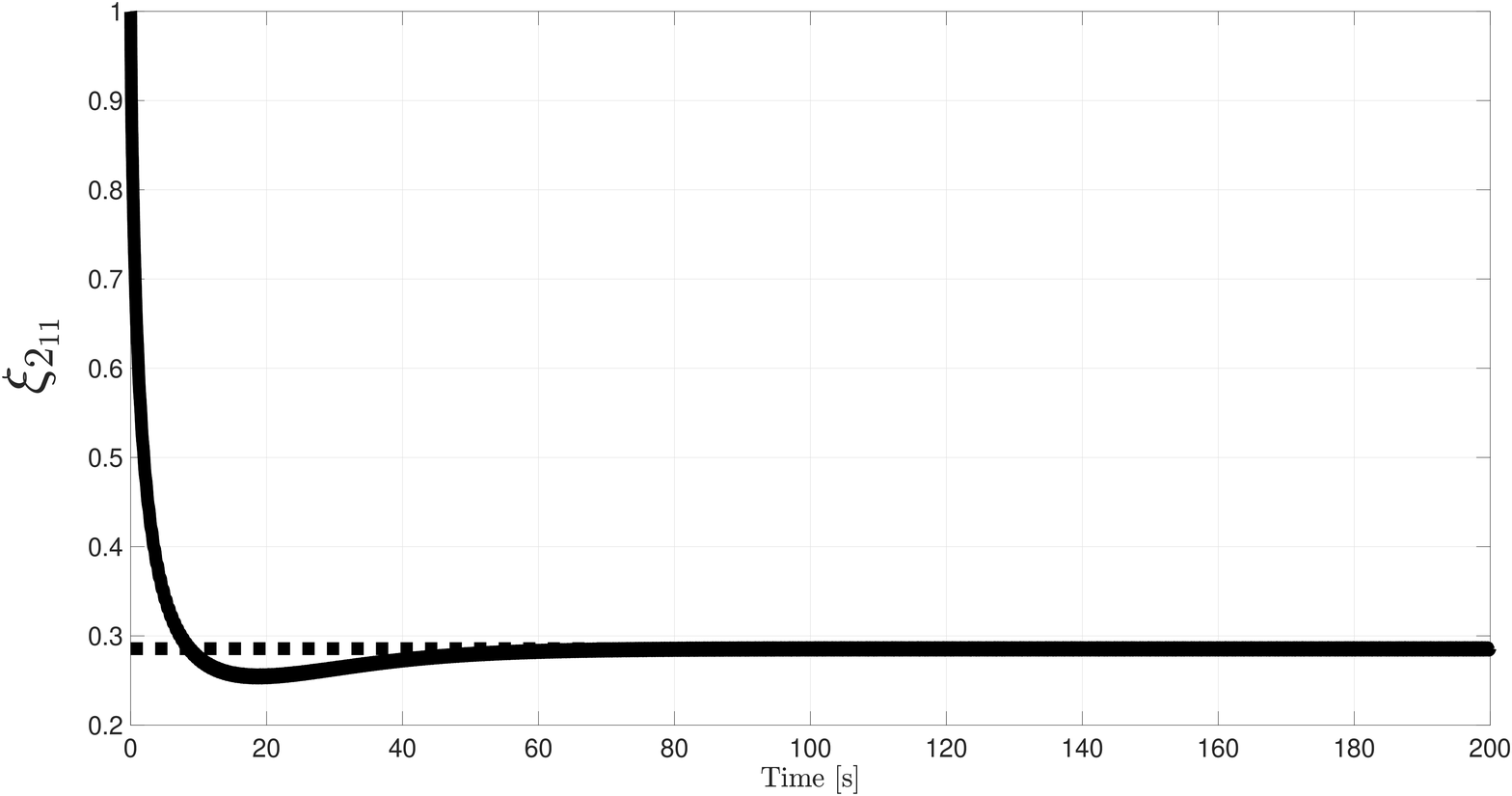} \includegraphics[width=0.24\textwidth]{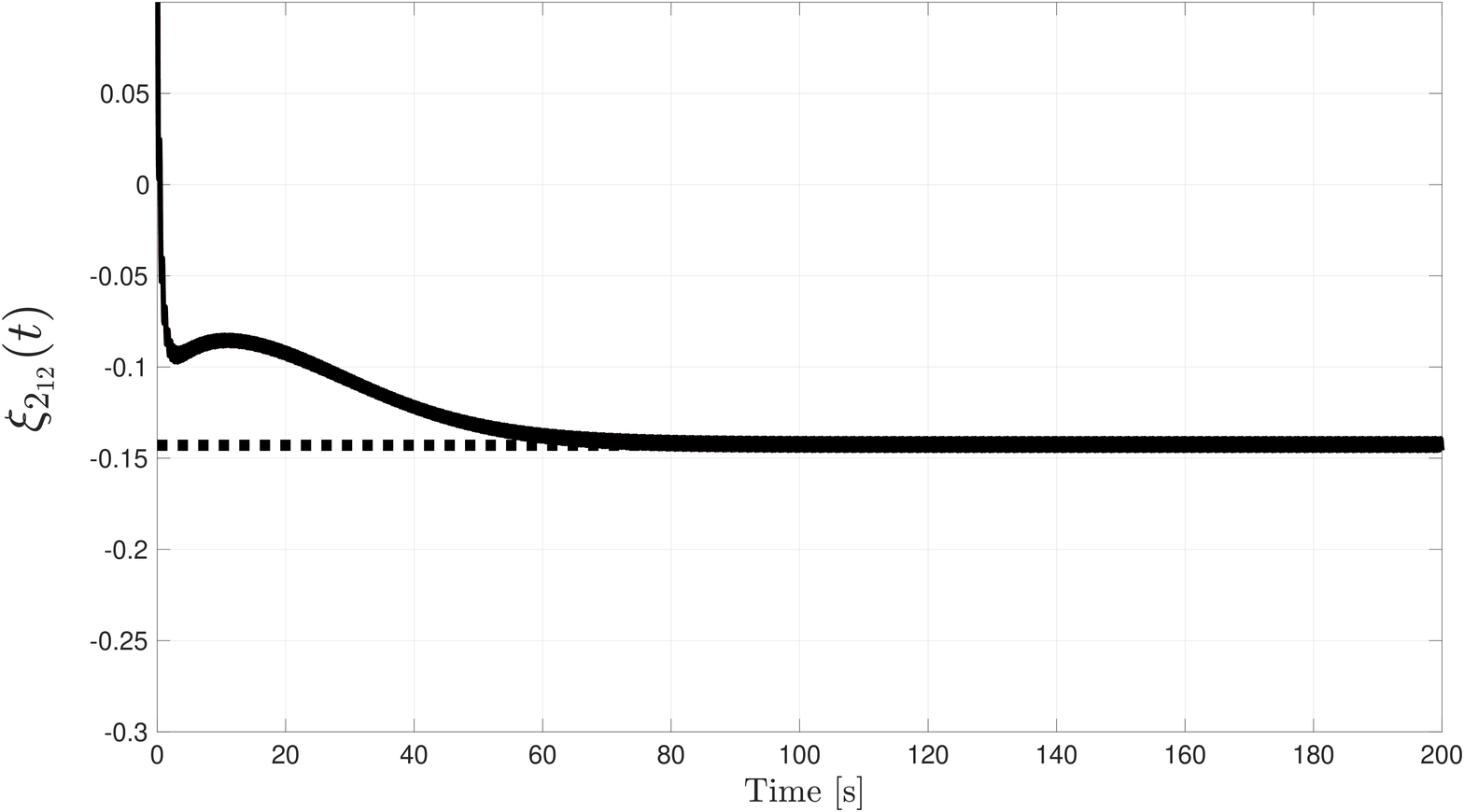}
  \caption{{Evolution in time of the component $H^{-1}_{11}$ of the inverse of the Hessian matrix.}}
  \label{Fig2}
 \end{centering}
\end{figure}

\begin{figure}[t!]
 \begin{centering}
  \includegraphics[width=0.24\textwidth]{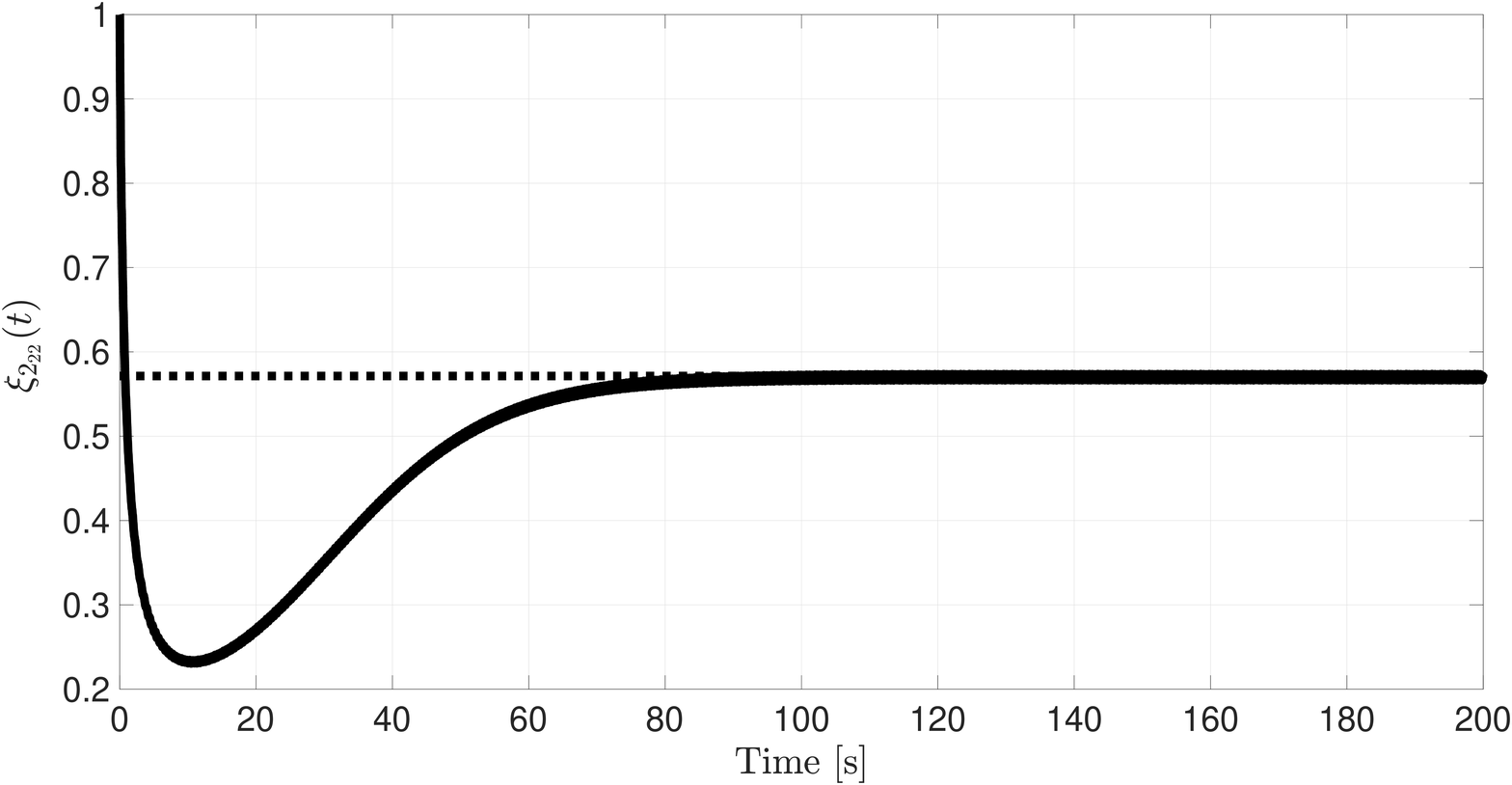} \includegraphics[width=0.24\textwidth]{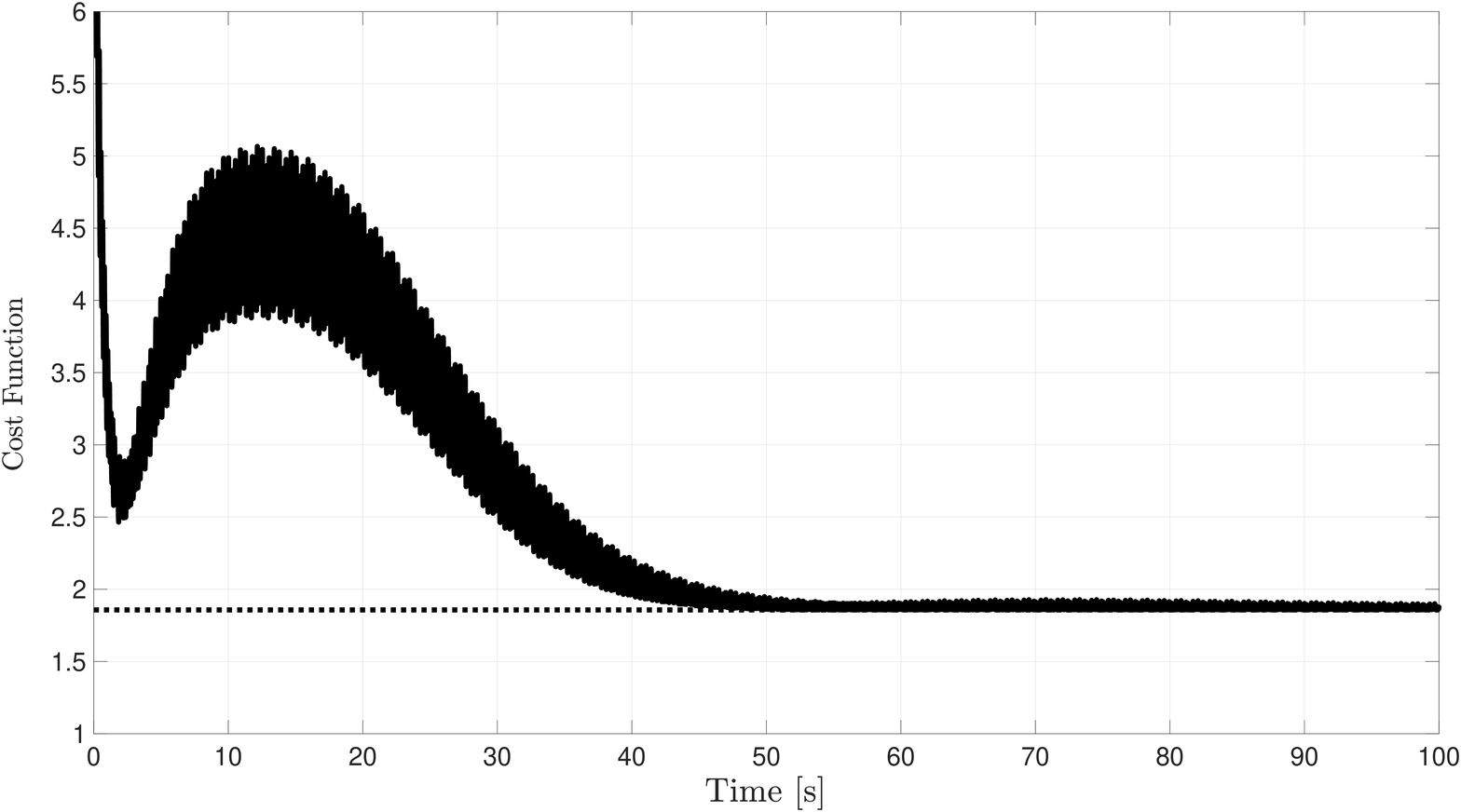}
  \caption{{Evolution in time of the component $H^{-1}_{22}$ of the inverse of the Hessian matrix.}}
  \label{Fig4}
 \end{centering}
\end{figure}

To illustrate the performance of the NFxTES, and to highlight the differences with respect to the the standard Newton-based extremum seeking controller of \cite{Newton}, we consider the quadratic function
\begin{equation}\label{cost_function}
\phi(z)=\frac{1}{2} z^\top H z+b^\top z+c,
\end{equation}
which satisfies $\nabla \phi(z)=Hz+b$ and $\nabla^2\phi(z)=H$. The parameters of \eqref{cost_function} are selected as
\begin{equation*}
H=\left[\begin{array}{cc}
4 & 1 \\ 1 & 2
\end{array}\right],~~~b=\left[\begin{array}{c}
-4\\
-6
\end{array}\right],~~c=11.
\end{equation*}
The inverse of the Hessian matrix is given by
\begin{equation}\label{hessiannumerical}
H^{-1}=\left[\begin{array}{cc}
0.2857 & -0.1429 \\ -0.1429 & 0.5714
\end{array}\right],
\end{equation}
and the function $\phi(z)$ has a global minimizer at the point
\begin{equation*}
z^*=-H^{-1}b=\left[\frac{2}{7},~\frac{20}{7}\right]^\top.
\end{equation*}
In order to find $z^*$ in fixed time, we implement the NFxTES dynamics with parameters $a=0.1$, $\varepsilon_1=0.1$ and $\varepsilon_2=10$. The constants $(k,q_1,q_2)$ were selected as $k=0.025$, $q_1=3$, and $q_2=1.5$, which generate an upper bound on the convergence time given by $T^*=123.4$.  We also simulate the Newton-based ES algorithm of \cite{Newton}, which has learning dynamics of the form $\dot{x}=-k\xi_1\xi_2$. To obtain a smooth approximation of $H^{-1}$, we have also implemented an additional low-pass filter that receives as input $\xi_1$ and $\xi_2$ and generates filtered outputs $\xi^f_1$ and $\xi^f_{2}$ that serve as inputs to the learning dynamics. As shown in \cite{Newton}, the incorporation of these filters does not affect the stability analysis of the algorithm. Figure \ref{Fig1} shows the trajectories generated by the NFxTES dynamics and by the classic Newton-based ES dynamics.  It can be observed that the NFxTES dynamics exhibit a much better transient performance in terms of less oscillations and faster convergence time to a neighborhood of $x^*$. On the other hand, Figures \ref{Fig2}-\ref{Fig4}, show the evolution of the components of the state $\xi_1$, which correspond to the entries of $H^{-1}$. As shown in the plots, these states converge to the true values of \eqref{hessiannumerical}.

\begin{figure}[t!]
 \begin{centering}
  \includegraphics[width=0.4\textwidth]{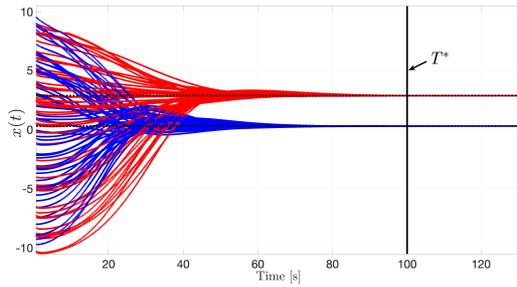}
  \caption{{Evolution in time of several solutions $x$ initialized in a neighborhood of the optimizers. The solid black line indicates the upper bound $T^*$ on the convergence time. }}
  \label{Fig6}
 \end{centering}
\end{figure}

 To further illustrate the fixed-time convergence property of the NFxTES dynamics, we have also simulated the case where the upper bound on the convergence time is fixed a priori as $T^*=100$, which can be obtained in the NFxTES dynamics by choosing $k=0.03085$, $q_1=3$, and $q_2=1.5$. Figure \ref{Fig6} shows the evolution in time of 50 different trajectories $x(t)$ initialized randomly in the set $[-10,10]\times[-10,10]$.  As it can be observed, the NFxTES dynamics guarantee convergence to a small neighborhood of $x^*$ before the prescribed time $T^*$. The simulations used the same parameters $(a,\varepsilon_1,\kappa_i)$ as in Figures 1-4, and $\varepsilon_2$ was selected as 6.25 to guarantee stability.

%
%\begin{figure}[t!]
% \begin{centering}
%  \includegraphics[width=0.48\textwidth]{cost.eps}
%  \caption{\small{Evolution in time of the cost function $\phi(z(t))$.}}
%  \label{Fig5}
% \end{centering}
%\end{figure}
%
\begin{figure}[t!]
 \begin{centering}
  \includegraphics[width=0.4\textwidth]{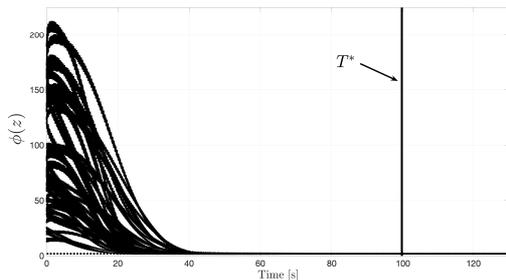}
  \caption{{Evolution in time of 50 trajectories of the cost function $\phi(z(t))$ generated by the 50 solutions $x$ shown in Figure \ref{Fig6}. The solid black line indicates the upper bound $T^*$ on the convergence time. }}
  \label{Fig7}
 \end{centering}
\end{figure}
%
%\begin{rem}
%It is important to note that in order to obtain then convergence result of Theorem \ref{main_theorem}, the parameters $(a,\varepsilon_1,\varepsilon_2)$ need to be appropriately tuned. In particular, these parameters depend on the constants $(q_1,q_2,c_1,c_2)$, which in turn fix the upper bound $T^*$ on the convergence time $T_c$. Thus, smaller values of $T^*$ may require smaller values for $(a,\varepsilon_1,\varepsilon_2)$. Since $\varepsilon_2$ is related to the frequency of the dither signal injected into the system, in order to implement in discrete time the NFxTES with a small $T^*$, one may need to use a sufficiently small step size to avoid aliasing issues. Therefore, there seems to be a tradeoff between achieving small finite convergence times and the computational complexity of the implemented discretized algorithm. \QEDB
%\end{rem}

%%%%%%%%%%%%%%%%%%%%%%%%%%%%%%%%
\section{Conclusions}
\label{sec_conclusions}
This paper presented a novel Newton-based extremum seeking controller that achieves fixed-time convergence in static maps, i.e., the convergence time is bounded by a constant that can be arbitrarily assigned by the designer. The learning dynamics of the extremum seeking algorithm implement a continuous vector field that receives as inputs the estimations of the gradient and the Hessian of the cost function $\phi(z)$, which are obtained by using only measurement of the cost. Local practical convergence in fixed-time was established by using tools from singular perturbation theory and averaging. The advantage of the method in comparison to the traditional Newton-based scheme was illustrated in numerical examples. %Future directions will study Newton-based schemes with fixed-time convergence results with arbitrarily large basins of attraction, as well as applications to source seeking in nonholonomic systems.

\bibliography{Biblio.bib}

\begin{thebibliography}{23}
\providecommand{\natexlab}[1]{#1}
\providecommand{\url}[1]{\texttt{#1}}
\providecommand{\urlprefix}{URL }
\expandafter\ifx\csname urlstyle\endcsname\relax
  \providecommand{\doi}[1]{doi:\discretionary{}{}{}#1}\else
  \providecommand{\doi}{doi:\discretionary{}{}{}\begingroup
  \urlstyle{rm}\Url}\fi

\bibitem[{Alvarez et~al.(2002)Alvarez, Attouch, Bolte, and
  Redont}]{NewtonSemiglobal}
Alvarez, F., Attouch, H., Bolte, J., and Redont, P. (2002).
\newblock A second-order gradient-like dissipative dynamical system with
  hessian-driven damping: {A}pplication to optimization and mechanics.
\newblock \emph{J. Math. Pures Appl.}, 81, 747--779.

\bibitem[{Andrieu et~al.(2008)Andrieu, Praly, and Astolfi}]{PralyFinite_Time}
Andrieu, V., Praly, L., and Astolfi, A. (2008).
\newblock Homogeneous approximation, recursive observer design, and output
  feedback.
\newblock \emph{SIAM J. Control Optim}, 47(4), 1814--1850.

\bibitem[{Ariyur and Krsti\'{c}(2003)}]{KrsticBookESC}
Ariyur, K.B. and Krsti\'{c}, M. (2003).
\newblock \emph{Real-Time Optimization by Extremum-Seeking Control}.
\newblock Wiley.

\bibitem[{Cruz-Zavala et~al.(2010)Cruz-Zavala, Moreno, and
  Fridman}]{JaimeFixed_Time}
Cruz-Zavala, E., Moreno, J.A., and Fridman, L. (2010).
\newblock Uniform second-order sliding mode observer for mechanical systems.
\newblock \emph{Proc. Int. Workshop Variable Struct. Syst}, 14--19.

\bibitem[{Engel and Kreisselmeier(2002)}]{finite_timeEngelTAC}
Engel, R. and Kreisselmeier, G. (2002).
\newblock A continuous-time observer which converges in finite time.
\newblock \emph{IEEE Transactions on Automatic and Control}, 47, 1202--1204.

\bibitem[{Garg and Panagou(2018)}]{fixed_time}
Garg, K. and Panagou, D. (2018).
\newblock Fixed-time stable gradient-flow schemes: Applications to
  continuous-time optimization.
\newblock \emph{arXiv:1808.10474}.

\bibitem[{Ghaffari et~al.(2012)Ghaffari, Krsti\'{c}, and
  Ne{\v{s}}i{\'c}}]{Newton}
Ghaffari, A., Krsti\'{c}, M., and Ne{\v{s}}i{\'c}, D. (2012).
\newblock Multivariable newton-based extremum seeking.
\newblock \emph{Automatica}, 48, 1759--1767.

\bibitem[{Ghaffari et~al.(2014)Ghaffari, Krsti\'{c}, and Seshagiri}]{PowerES}
Ghaffari, A., Krsti\'{c}, M., and Seshagiri, S. (2014).
\newblock Power optimization for photovoltaic micro-converters using
  multivariable gradient-based extremum-seeking.
\newblock \emph{IEEE Transactions on Control System Technology}, 22(6),
  2141--2149.

\bibitem[{Grushkovskaya et~al.(2017)Grushkovskaya, Durr, Ebenbauer, and
  Zuyev}]{Grushkovskaya2017}
Grushkovskaya, V., Durr, H., Ebenbauer, C., and Zuyev, A. (2017).
\newblock Extremum seeking for time-varying functions using lie bracket
  approximations.
\newblock \emph{20th {W}orld {C}ongress}, 50(1), 522--5528.

\bibitem[{Guay and Zhang(2003)}]{Guay:03}
Guay, M. and Zhang, T. (2003).
\newblock Adaptive extremum seeking control of nonlinear dynamic systems with
  parametric uncertainties.
\newblock \emph{Automatica}, 39, 1283--1293.

\bibitem[{Li et~al.(2017)Li, Yu, Zhou, and Rein}]{WeiRenFixed}
Li, C., Yu, X., Zhou, X., and Rein, W. (2017).
\newblock A fixed time distributed optimization: A sliding mode perspective.
\newblock \emph{In Proc. of 43rd Annual Conference of the IEEE Industrial
  Electronics Society}, 8201--8207.

\bibitem[{Ne\u{s}i\'{c} et~al.(2010)Ne\u{s}i\'{c}, Tan, Moase, and
  Manzie}]{DerivativesESC}
Ne\u{s}i\'{c}, D., Tan, Y., Moase, W.H., and Manzie, C. (2010).
\newblock A unifying approach to extremum seeking: Adaptive schemes based on
  estimation of derivatives.
\newblock \emph{49th IEEE Conference on Decision and Control}, 4625--4630.

\bibitem[{Oliveira et~al.(2017)Oliveira, Krsti\'{c}, and
  Tsubakino}]{ThiagoKrstic}
Oliveira, T., Krsti\'{c}, M., and Tsubakino, D. (2017).
\newblock Extremum seeking for static maps with delays.
\newblock \emph{IEEE Trans. Autom. Control}, 62(4), 1911--1926.

\bibitem[{Polyakov(2012)}]{Fixed_timeTAC}
Polyakov, A. (2012).
\newblock Nonlinear feedback design for fixed-time stabilization of linear
  control systems.
\newblock \emph{IEEE Transactions on Automatic and Control}, 57(8), 2106--2110.

\bibitem[{Poveda and Krsti\'{c}(2020{\natexlab{a}})}]{PovedaKrsticACC20}
Poveda, J.I. and Krsti\'{c}, M. (2020{\natexlab{a}}).
\newblock Gradient-based fixed-time extremum seeking.
\newblock \emph{In Proc. of American Control Conference.}, 2838--2843.

\bibitem[{Poveda and Krsti\'{c}(2020{\natexlab{b}})}]{fixed_timeES_arxiv}
Poveda, J.I. and Krsti\'{c}, M. (2020{\natexlab{b}}).
\newblock Non-smooth extremum seeking control with user-prescribed convergence.
\newblock \emph{IEEE Transactions on Automatic and Control, provisionally
  accepted}.

\bibitem[{Poveda and Li(2021)}]{zero_order_poveda_Lina}
Poveda, J.I. and Li, N. (2021).
\newblock Robust hybrid zero-order optimization algorithms with acceleration
  via averaging in time.
\newblock \emph{Automatica}, 123.

\bibitem[{Poveda and Teel(2017)}]{Poveda:16}
Poveda, J.I. and Teel, A.R. (2017).
\newblock A framework for a class of hybrid extremum seeking controllers with
  dynamic inclusions.
\newblock \emph{Automatica}, 76, 113--126.

\bibitem[{Rios and Teel(2018)}]{RiosTeel18}
Rios, H. and Teel, A.R. (2018).
\newblock A hybrid fixed-time observer for state estimation of linear systems.
\newblock \emph{Automatica}, 87, 103--112.

\bibitem[{Romero and Benosman(2019)}]{romeronips}
Romero, O. and Benosman, M. (2019).
\newblock \emph{Finite-Time Convergence of Continuous-Time Optimization
  Algorithms via Differential Inclusions}.
\newblock Neurips, to appear.

\bibitem[{Suttner(2019)}]{ESC_Suttner}
Suttner, R. (2019).
\newblock Extremum seeking control with an adaptive dither signal.
\newblock \emph{Automatica}, 101, 214--222.

\bibitem[{Tan et~al.(2006)Tan, Ne\v{s}i\'{c}, and Mareels}]{tan06Auto}
Tan, Y., Ne\v{s}i\'{c}, D., and Mareels, I.M. (2006).
\newblock On non-local stability properties of extremum seeking control.
\newblock \emph{Automatica}, 42(6), 889--903.

\bibitem[{Teel and Popovic(2001)}]{Teel:01}
Teel, A.R. and Popovic, D. (2001).
\newblock Solving smooth and nonsmooth multivariable extremum seeking problems
  by the methods of nonlinear programming.
\newblock \emph{In proc. of American Control Conference}, 2394--2399.

\end{thebibliography}
\end{document}